\documentclass[12pt]{article}%
\usepackage{graphicx}
\usepackage[intlimits]{amsmath}
\usepackage{latexsym}
\usepackage{amsfonts}
\usepackage{amssymb}%
 \makeatletter
\newfont{\footsc}{cmcsc10 at 8truept}
\newfont{\footbf}{cmbx10 at 8truept}
\newfont{\footrm}{cmr10 at 10truept}
\newtheorem{theorem}{Theorem}

\newtheorem{fact}[theorem]{Fact}

\newenvironment{proof}[1][Proof]{\noindent{\textbf {#1}  }}  {\hfill$\Box$\bigskip}

\begin{document}

\title{Spectral saturation: inverting the spectral Tur\'{a}n theorem}
\author{Vladimir Nikiforov\\{\small Department of Mathematical Sciences, University of Memphis, Memphis TN
38152}\\{\small email: vnikifrv@memphis.edu}}
\maketitle

\begin{abstract}
Let $\mu\left(  G\right)  $ be the largest eigenvalue of a graph $G$ and
$T_{r}\left(  n\right)  $ be the $r$-partite Tur\'{a}n graph of order $n.$

We prove that if $G$ is a graph of order $n$ with $\mu\left(  G\right)
>\mu\left(  T_{r}\left(  n\right)  \right)  ,$ then $G$ contains varoius large
supergraphs of the complete graph of order $r+1,$ e.g., the complete
$r$-partite graph with all parts of size $\log n$ with an edge added to the
first part.

We also give corresponding stability results.\medskip


\end{abstract}

\section{Introduction}

This note is part of an ongoing project aiming to build extremal graph theory
on spectral grounds, see, e.g., \cite{BoNi07}, \cite{Nik02,Nik07h}.

Let $\mu\left(  G\right)  $ be the largest adjacency eigenvalue of a graph $G$
and $T_{r}\left(  n\right)  $ be the $r$-partite Tur\'{a}n graph of order $n.$
The spectral Tur\'{a}n theorem \cite{Nik07d} implies that if $G$ is a graph of
order $n$ with $\mu\left(  G\right)  >\mu\left(  T_{r}\left(  n\right)
\right)  ,$ then $G$ contains a $K_{r+1}$, the complete graph of order $r+1.$

On the other hand, it is known (e.g., \cite{BoNi04}, \cite{Erd63},
\cite{ErSi73}, \cite{Nik07b}) that if $e\left(  G\right)  >e\left(
T_{r}\left(  n\right)  \right)  ,$ then $G$ contains large supergraphs of
$K_{r+1}.$

It turns out that essentially the same results also follow from $\mu\left(
G\right)  >\mu\left(  T_{r}\left(  n\right)  \right)  .$

Recall first a family of graphs, studied initially by Erd\H{o}s \cite{Erd69}
and recently in \cite{BoNi04}: an\emph{ }$r$\emph{-joint }of size $t$ is the
union of $t$ distinct $r$-cliques sharing an edge. Write $js_{r}\left(
G\right)  $ for the maximum size of an $r$-joint in a graph $G.$ Erd\H{o}s
\cite{Erd69}, Theorem 3', showed that:\medskip

\emph{If }$G$\emph{ is a graph of sufficiently large order }$n$\emph{
satisfies }$e\left(  G\right)  >e\left(  T_{r}\left(  n\right)  \right)
$\emph{, then} $js_{r+1}\left(  G\right)  >n^{r-1}/\left(  10\left(
r+1\right)  \right)  ^{6\left(  r+1\right)  }.$\medskip

Here is a explicit spectral analogue of this result.

\begin{theorem}
\label{th1}Let $r\geq2,$ $n>r^{15},$ and $G$ be a graph of order $n.$ If
$\mu\left(  G\right)  >\mu\left(  T_{r}\left(  n\right)  \right)  ,$ then
$js_{r+1}\left(  G\right)  >n^{r-1}/r^{2r+4}.$
\end{theorem}

Erd\H{o}s \cite{Erd63} introduced yet another graph related to Tur\'{a}n's
theorem: let $K_{r}^{+}\left(  s_{1},\ldots,s_{r}\right)  $ be the complete
$r$-partite graph with parts of size $s_{1}\geq2,s_{2},\ldots,s_{r},$ with an
edge added to the first part. The extremal results about this graph given in
\cite{Erd63} and \cite{ErSi73} were recently extended in \cite{Nik07b}
to:\medskip

\emph{Let }$r\geq2,$\emph{ }$2/\ln n\leq c\leq r^{-\left(  r+7\right)  \left(
r+1\right)  },$\emph{ and }$G$\emph{ be a graph of order }$n$\emph{. If }%
$G$\emph{ has }$t_{r}\left(  n\right)  +1$\emph{ edges, then }$G$\emph{
contains a }$K_{r}^{+}\left(  \left\lfloor c\ln n\right\rfloor ,\ldots
,\left\lfloor c\ln n\right\rfloor ,\left\lceil n^{1-\sqrt{c}}\right\rceil
\right)  .$\medskip

Here we give a similar spectral extremal result.

\begin{theorem}
\label{th2} Let $r\geq2,$ $2/\ln n\leq c\leq r^{-\left(  2r+9\right)  \left(
r+1\right)  },$ and $G$ be a graph of order $n.$ If $\mu\left(  G\right)
>\mu\left(  T_{r}\left(  n\right)  \right)  ,$ then $G$ contains a $K_{r}%
^{+}\left(  \left\lfloor c\ln n\right\rfloor ,\ldots,\left\lfloor c\ln
n\right\rfloor ,\left\lceil n^{1-\sqrt{c}}\right\rceil \right)  .$
\end{theorem}

As an easy consequence of Theorem \ref{th2} we obtain

\begin{theorem}
\label{th3} Let $r\geq2,$ $c=r^{-\left(  2r+9\right)  \left(  r+1\right)  },$
$n\geq e^{2/c},$ and $G$ be a graph of order $n.$ If $\mu\left(  G\right)
>\mu\left(  T_{r}\left(  n\right)  \right)  ,$ then $G$ contains a $K_{r}%
^{+}\left(  \left\lfloor c\ln n\right\rfloor ,\ldots,\left\lfloor c\ln
n\right\rfloor \right)  .$
\end{theorem}

Theorems \ref{th1}, \ref{th2}, and \ref{th3} have corresponding stability results.

\begin{theorem}
\label{th1.2}Let $r\geq2,$ $0<b<2^{-10}r^{-6},$ $n\geq r^{20},$ and $G$ be a
graph of order $n.$ If $\mu\left(  G\right)  >\left(  1-1/r-b\right)  n,$ then
$G$ satsisfies one of the conditions:

(a) $js_{r+1}\left(  G\right)  >n^{r-1}/r^{2r+5};$

(b) $G$ contains an induced $r$-partite subgraph $G_{0}$ of order at least
$\left(  1-4b^{1/3}\right)  n$ with minimum degree $\delta\left(
G_{0}\right)  >\left(  1-1/r-7b^{1/3}\right)  n.${}
\end{theorem}

\begin{theorem}
\label{th2.2}Let $r\geq2,$ $2/\ln n\leq c\leq r^{-\left(  2r+9\right)  \left(
r+1\right)  }/2,$ $0<b<2^{-10}r^{-6},$ and $G$ be a graph of order $n.$ If
$\mu\left(  G\right)  >\left(  1-1/r-b\right)  n$, then $G$ satsisfies one of
the conditions:

(a) $G$ contains a $K_{r}^{+}\left(  \left\lfloor c\ln n\right\rfloor
,\ldots,\left\lfloor c\ln n\right\rfloor ,\left\lceil n^{1-2\sqrt{c}%
}\right\rceil \right)  ;$

(b) $G$ contains an induced $r$-partite subgraph $G_{0}$ of order at least
$\left(  1-4b^{1/3}\right)  n$ with minimum degree $\delta\left(
G_{0}\right)  >\left(  1-1/r-7b^{1/3}\right)  n.$
\end{theorem}

\begin{theorem}
\label{th3.2}Let $r\geq2,$ $c=r^{-\left(  2r+9\right)  \left(  r+1\right)
}/2,$ $0<b<2^{-10}r^{-6},$ $n\geq e^{2/c},$ and $G$ be a graph of order $n.$
If $\mu\left(  G\right)  >\left(  1-1/r-b\right)  n$, then one of the
following conditions holds:

(a) $G$ contains a $K_{r}^{+}\left(  \left\lfloor c\ln n\right\rfloor
,\ldots,\left\lfloor c\ln n\right\rfloor \right)  ;$

(b) $G$ contains an induced $r$-partite subgraph $G_{0}$ of order at least
$\left(  1-4b^{1/3}\right)  n$ with minimum degree $\delta\left(
G_{0}\right)  >\left(  1-1/r-7b^{1/3}\right)  n.$
\end{theorem}

\subsubsection*{Remarks}

\begin{itemize}
\item[-] Obviously Theorems \ref{th1}, \ref{th2}, and \ref{th3} are tight
since $T_{r}\left(  n\right)  $ contains no $\left(  r+1\right)  $-cliques.

\item[-] Theorems \ref{th2}, \ref{th3}, \ref{th2.2}, and \ref{th3.2} are
essentially best possible since for every $\varepsilon>0,$ choosing randomly a
graph $G$ of order $n$ with $e\left(  G\right)  =\left\lceil \left(
1-\varepsilon\right)  n^{2}/2\right\rceil $ edges we see that $\mu\left(
G\right)  >\left(  1-\varepsilon\right)  n,$ but $G$ contains no $K_{2}\left(
c\ln n,c\ln n\right)  $ for some $c>0,$ independent of $n.$

\item[-] Theorem \ref{th1} implies in turn spectral versions of other known
results, like Theorem 3.8 in \cite{EFS94}:

\emph{Every graph }$G$\emph{ of order }$n$\emph{ with }$\mu\left(  G\right)
>\mu\left(  T_{r}\left(  n\right)  \right)  $\emph{ contains }$cn$\emph{
distinct }$\left(  r+1\right)  $\emph{-cliques sharing an }$r$\emph{-clique,
where }$c>0$\emph{ is independent of }$n.$

\item[-] The relations between $c$ and $n$ in Theorems \ref{th2} and
\ref{th2.2} need explanation. First, for fixed $c,$ they show how large must
be $n$ to get valid conclusions. But, in fact, the relations are subtler, for
$c$ itself may depend on $n,$ e.g., letting $c=1/\ln\ln n,$ the conclusions
are meaningful for sufficiently large $n.$

\item[-] Note that, in Theorems \ref{th2} and \ref{th2.2}, if the conclusion
holds for some $c,$ it holds also for $0<c^{\prime}<c,$ provided $n$ is
sufficiently large;

\item[-] The stability conditions \emph{(b) }in Theorems \ref{th1.2},
\ref{th2.2}, and \ref{th3.2} are stronger than the conditions in the stability
theorems of \cite{Erd68}, \cite{Sim68} and \cite{Nik07a}. Indeed, in all these
theorems, condition \emph{(ii)} implies that $G_{0}$ is an induced, almost
balanced, and almost complete $r$-partite graph containing almost all the
vertices of $G;$

\item[-] The exponents $1-\sqrt{c}$ and $1-2\sqrt{c}$ in Theorems \ref{th2}
and \ref{th2.2} are far from the best ones, but are simple.
\end{itemize}

The next section contains notation and results needed to prove the theorems.
The proofs are presented in Section \ref{pf}.

\section{\label{pr}Preliminary results}

Our notation follows \cite{Bol98}. Given a graph $G,$ we write:\medskip

- $V\left(  G\right)  $ for the vertex set of $G$ and $\left\vert G\right\vert
$ for $\left\vert V\left(  G\right)  \right\vert ;$

- $E\left(  G\right)  $ for the edge set of $G$ and $e\left(  G\right)  $ for
$\left\vert E\left(  G\right)  \right\vert ;$

- $d\left(  u\right)  $ for the degree of a vertex $u;$

- $\delta\left(  G\right)  $ for the minimum degree of $G;$

- $k_{r}\left(  G\right)  $ for the number of $r$-cliques of $G;$

- $K_{r}\left(  s_{1},\ldots,s_{r}\right)  $ for the complete $r$-partite
graph with parts of size $s_{1},\ldots,s_{r}.$\bigskip

The following facts play crucial roles in our proofs.

\begin{fact}
[\cite{Nik07d}, Theorem 1]\label{STT}Every graph $G$ of order $n$ with
$\mu\left(  G\right)  >\mu\left(  T_{r}\left(  n\right)  \right)  $ contains a
$K_{r+1}.\hfill\square$
\end{fact}

\begin{fact}
[\textbf{\cite{Nik07c}, Theorem 5}]\label{thv3}Let $\ 0<\alpha\leq1/4,$
\ $0<\beta\leq1/2,$ \ $1/2-\alpha/4\leq\gamma<1,$ $\ K\geq0,$ $\ n\geq\left(
42K+4\right)  /\alpha^{2}\beta$, and $G$ be a graph of order $n.$ If
\[
\mu\left(  G\right)  >\gamma n-K/n\text{ \ \ and \ \ }\delta\left(  G\right)
\leq\left(  \gamma-\alpha\right)  n,
\]
then $G$ contains an induced subgraph $H$ satisfying$\ \left\vert H\right\vert
\geq\left(  1-\beta\right)  n$ and one of the conditions:

(a) $\mu\left(  H\right)  >\gamma\left(  1+\beta\alpha/2\right)  \left\vert
H\right\vert ;$

(b) $\mu\left(  H\right)  >\gamma\left\vert H\right\vert $ and $\delta\left(
H\right)  >\left(  \gamma-\alpha\right)  \left\vert H\right\vert
.\hfill\square$
\end{fact}

\begin{fact}
[\cite{BoNi04}, Lemma 6]\label{leKd} Let $r\geq2$ and $G$ be graph a of order
$n.$ If $G$ contains a $K_{r+1}$ and $\ \delta\left(  G\right)  >\left(
1-1/r-1/r^{4}\right)  n,$ then $js_{r+1}\left(  G\right)  >n^{r-1}%
/r^{r+3}.\hfill\square$
\end{fact}

\begin{fact}
[\cite{BoNi07}, Theorem 2]\label{leNSMM}If $r\geq2$ and $G$ is a graph of
order $n,$ then%
\[
k_{r}\left(  G\right)  \geq\left(  \frac{\mu\left(  G\right)  }{n}-1+\frac
{1}{r}\right)  \frac{r\left(  r-1\right)  }{r+1}\left(  \frac{n}{r}\right)
^{r+1}.
\]
$\hfill\square$
\end{fact}

\begin{fact}
[\cite{BoNi07}, Theorem 4]\label{tstab} Let $r\geq2,$ $0\leq b\leq
2^{-10}r^{-6},$ and $G$ be a graph of order $n.$ If $G$ contains no $K_{r+1}$
and $\mu\left(  G\right)  \geq\left(  1-1/r-b\right)  n,$ then $G$ contains an
induced $r$-partite graph $G_{0}$ satisfying $\left\vert G_{0}\right\vert
\geq\left(  1-3c^{1/3}\right)  n$ and $\delta\left(  G_{0}\right)  >\left(
1-1/r-6c^{1/3}\right)  n.\hfill\square$
\end{fact}

\begin{fact}
[\cite{Nik07b}, Theorem 6]\label{thv4}Let $r\geq2,$ $2/\ln n\leq c\leq
r^{-\left(  r+8\right)  r},$ and $g$ is a graph $G$ of order $n.$ If $G$
contains a $K_{r+1}$ and $\delta\left(  G\right)  >\left(  1-1/r-1/r^{4}%
\right)  n,$ then $G$ contains a $K_{r}^{+}\left(  \left\lfloor c\ln
n\right\rfloor ,\ldots,\left\lfloor c\ln n\right\rfloor ,\left\lceil
n^{1-cr^{3}}\right\rceil \right)  .\hfill\square$
\end{fact}

\begin{fact}
[\cite{Nik07}, Theorem 1]\label{ES}Let $r\geq2,$ $c^{r}\ln n\geq1,$ and $G$ be
a graph of order $n$. If $k_{r}\left(  G\right)  \geq cn^{r},$ then $G$
contains a $K_{r}\left(  s,\ldots s,t\right)  $ with $s=\left\lfloor c^{r}\ln
n\right\rfloor $ and $t>n^{1-c^{r-1}}.\hfill\square$
\end{fact}

\begin{fact}
\label{tsize}The number of edges of $T_{r}\left(  n\right)  $ satisfies
$2e\left(  T_{r}\left(  n\right)  \right)  \geq\left(  1-1/r\right)
n^{2}-r/4.\hfill\square$
\end{fact}

\section{\label{pf}Proofs}

Below we prove Theorems \ref{th1}, \ref{th2}, \ref{th1.2}, and \ref{th2.2}. We
omit the proofs of Theorems \ref{th3} and \ref{th3.2} since they are easy
consequences of Theorems \ref{th2} and \ref{th2.2}.

All proofs have similar simple structure and follow from the facts listed above.

\subsubsection*{Proof of Theorem \ref{th1}}

\begin{proof}
[ ]Let $G$ be a graph of order $n$ with $\mu\left(  G\right)  >\mu\left(
T_{r}\left(  n\right)  \right)  ;$ thus, by Fact \ref{STT}, $G$ contains a
$K_{r+1}.$ If%
\begin{equation}
\delta\left(  G\right)  >\left(  1-r^{-1}-r^{-4}\right)  n, \label{mind}%
\end{equation}
then, by Fact \ref{leKd}, $js_{r+1}\left(  G\right)  >n^{r-1}/r^{r+3},$
completing the proof.

Thus, we shall assume that (\ref{mind}) fails. Then, letting
\begin{equation}
\alpha=1/r^{4},\text{ \ \ }\beta=1/2,\text{ \ \ }\gamma=1-1/r,\text{
\ \ }K=r/4, \label{co1}%
\end{equation}
we see that
\begin{equation}
\delta\left(  G\right)  \leq\left(  \gamma-\alpha\right)  n \label{co2}%
\end{equation}
and also, in view of Fact \ref{tsize},
\begin{equation}
\mu\left(  G\right)  >\mu\left(  T_{r}\left(  n\right)  \right)  \geq2e\left(
T_{r}\left(  n\right)  \right)  /n\geq\left(  1-1/r\right)  n-r/4n=\gamma
n-K/n. \label{co3}%
\end{equation}
Given (\ref{co1}), (\ref{co2}) and (\ref{co3}), Theorem \ref{thv3} implies
that, for $n\geq r^{15}$, $G$ contains an induced subgraph $H$ satisfying
$\left\vert H\right\vert \geq n/2$ and one of the conditions:

\emph{(i)} $\mu\left(  H\right)  >\left(  1-1/r+1/\left(  4r^{4}\right)
\right)  \left\vert H\right\vert ;$

\emph{(ii)} $\mu\left(  H\right)  >\left(  1-1/r\right)  \left\vert
H\right\vert $ and $\delta\left(  H\right)  >\left(  1-1/r-1/r^{4}\right)
\left\vert H\right\vert .$

If condition \emph{(i)} holds, Fact \ref{leNSMM} gives
\[
k_{r+1}\left(  H\right)  >\left(  \frac{\mu\left(  H\right)  }{\left\vert
H\right\vert }-1-\frac{1}{r}\right)  \frac{r\left(  r-1\right)  }{r+1}\left(
\frac{\left\vert H\right\vert }{r}\right)  ^{r+1}>\frac{r\left(  r-1\right)
}{4r^{4}\left(  r+1\right)  }\left(  \frac{\left\vert H\right\vert }%
{r}\right)  ^{r+1},
\]
and so,
\begin{align*}
js_{r+1}\left(  G\right)   &  \geq js_{r+1}\left(  H\right)  \geq\binom
{r+1}{2}\frac{k_{r+1}\left(  H\right)  }{e\left(  H\right)  }>r\left(
r+1\right)  \frac{k_{r+1}\left(  H\right)  }{\left\vert H\right\vert ^{2}}\\
&  >\frac{r\left(  r+1\right)  r\left(  r-1\right)  }{4r^{4}\left(
r+1\right)  r^{r+1}}\left\vert H\right\vert ^{r-1}>\frac{1}{4r^{r+3}%
}\left\vert H\right\vert ^{r-1}\geq\frac{1}{2^{r+1}r^{r+3}}n^{r-1}>\frac
{1}{r^{2r+4}}n^{r-1},
\end{align*}
completing the proof.

If condition \emph{(ii)} holds, then $H$ contains a $K_{r+1}$; thus, by Fact
\ref{leKd}, $js_{r+1}\left(  H\right)  >\left\vert H\right\vert ^{r-1}%
/r^{r+3}.$ To complete the proof, notice that
\[
js_{r+1}\left(  G\right)  >js_{r+1}\left(  H\right)  >\frac{\left\vert
H\right\vert ^{r-1}}{r^{r+3}}\geq\frac{1}{2^{r-1}r^{r+3}}n^{r-1}>\frac
{1}{r^{2r+4}}n^{r-1}.
\]

\end{proof}

\bigskip

\subsubsection*{Proof of Theorem \ref{th2}}

\begin{proof}
[ ]Let $G$ be a graph of order $n$ with $\mu\left(  G\right)  >\mu\left(
T_{r}\left(  n\right)  \right)  ;$ thus, by Fact \ref{STT}, $G$ contains a
$K_{r+1}.$ If%
\begin{equation}
\delta\left(  G\right)  >\left(  1-1/r-1/r^{4}\right)  n, \label{mind1}%
\end{equation}
then, by Fact \ref{thv4}, $G$ contains a $K_{r}^{+}\left(  \left\lfloor c\ln
n\right\rfloor ,\ldots,\left\lfloor c\ln n\right\rfloor ,\left\lceil
n^{1-cr^{3}}\right\rceil \right)  ,$ completing the proof, in view of
$cr^{3}<\sqrt{c}.$

Thus, we shall assume that (\ref{mind1}) fails. Then, letting
\begin{equation}
\alpha=1/r^{4},\text{ \ \ }\beta=1/2,\text{ \ \ }\gamma=1-1/r,\text{
\ \ }K=r/4, \label{co4}%
\end{equation}
we see that
\begin{equation}
\delta\left(  G\right)  \leq\left(  \gamma-\alpha\right)  n \label{co5}%
\end{equation}
and also, in view of Fact \ref{tsize},
\begin{equation}
\mu\left(  G\right)  >\mu\left(  T_{r}\left(  n\right)  \right)  \geq2e\left(
T_{r}\left(  n\right)  \right)  /n\geq\left(  1-1/r\right)  n-r/4n=\gamma
n-K/n. \label{co6}%
\end{equation}
Given (\ref{co4}), (\ref{co5}) and (\ref{co6}), Theorem \ref{thv3} implies
that, for $n>r^{15}$, $G$ contains an induced subgraph $H$ satisfying
$\left\vert H\right\vert \geq n/2$ and one of the conditions:

\emph{(i)} $\mu\left(  H\right)  >\left(  1-1/r+1/\left(  4r^{4}\right)
\right)  \left\vert H\right\vert ;$

\emph{(ii)} $\mu\left(  H\right)  >\left(  1-1/r\right)  \left\vert
H\right\vert $ and $\delta\left(  H\right)  >\left(  1-1/r-1/r^{4}\right)
\left\vert H\right\vert .$

If condition \emph{(i)} holds, Fact \ref{leNSMM} gives
\begin{align*}
k_{r+1}\left(  H\right)   &  >\left(  \frac{\mu\left(  H\right)  }{\left\vert
H\right\vert }-1-\frac{1}{r}\right)  \frac{r\left(  r-1\right)  }{r+1}\left(
\frac{\left\vert H\right\vert }{r}\right)  ^{r+1}>\frac{r\left(  r-1\right)
}{4r^{4}\left(  r+1\right)  }\left(  \frac{\left\vert H\right\vert }%
{r}\right)  ^{r+1}\\
&  >\frac{1}{2^{r+3}r^{r+4}\left(  r+1\right)  }n^{r+1}>\frac{1}{r^{2r+9}%
}n^{r+1}\geq c^{1/\left(  r+1\right)  }n^{r+1}.
\end{align*}
Thus, by Fact \ref{ES}, $G$ contains a $K_{r+1}\left(  s,\ldots,s,t\right)  $
with $s=\left\lfloor c\ln n\right\rfloor $ and $t>n^{1-c^{r/\left(
r+1\right)  }}>n^{1-\sqrt{c}}.$ Then, obviously, $G$ contains a $K_{r}%
^{+}\left(  \left\lfloor c\ln n\right\rfloor ,\ldots,\left\lfloor c\ln
n\right\rfloor ,\left\lceil n^{1-\sqrt{c}}\right\rceil \right)  ,$ completing
the proof.

If condition \emph{(ii)} holds, then $H$ contains a $K_{r+1}$; thus, by Fact
\ref{thv4}, $H$ contains a
\[
K_{r}^{+}\left(  \left\lfloor 2c\ln\left\vert H\right\vert \right\rfloor
,\ldots,\left\lfloor 2c\ln\left\vert H\right\vert \right\rfloor ,\left\lceil
\left\vert H\right\vert ^{1-2cr^{3}}\right\rceil \right)  .
\]
To complete the proof, note that $2c\ln\left\vert H\right\vert \geq2c\ln
\frac{n}{2}>c\ln n$ and%
\[
\left\vert H\right\vert ^{1-2cr^{3}}\geq\left(  \frac{n}{2}\right)
^{1-2cr^{3}}\geq\frac{1}{2}n^{1-2cr^{3}}>n^{1-\sqrt{c}}.
\]

\end{proof}

\bigskip

\begin{proof}
[\textbf{Proof of Theorem \ref{th1.2}}]Let $G$ be a graph of order $n$ with
$\mu\left(  G\right)  >\left(  1-1/r-b\right)  n$. If $G$ contains no
$K_{r+1},$ then condition \emph{(b) }follows from Fact \ref{tstab}; thus we
assume that $G$ contains a $K_{r+1}.$ If
\begin{equation}
\delta\left(  G\right)  >\left(  1-1/r-1/r^{4}\right)  n, \label{mind2}%
\end{equation}
then Fact \ref{leKd} implies condition \emph{(a)}.

Thus, we shall assume that (\ref{mind2}) fails. Then, letting
\begin{equation}
\alpha=1/r^{4}-b,\text{ \ \ }\beta=4b/\alpha,\text{ \ \ }\gamma=1-1/r-b,\text{
\ \ }K=0, \label{co7}%
\end{equation}
we easily see that
\begin{equation}
\beta=\frac{4b}{1/r^{4}-b}\leq\frac{1}{2},\text{ \ \ \ }\delta\left(
G\right)  \leq\left(  \gamma-\alpha\right)  n, \label{co8}%
\end{equation}
and
\begin{equation}
\mu\left(  G\right)  >\left(  1-1/r-b\right)  n=\gamma n. \label{co9}%
\end{equation}
Given (\ref{co7}), (\ref{co8}) and (\ref{co9}), Theorem \ref{thv3} implies
that, for $n\geq r^{20}$, $G$ contains an induced subgraph $H$ satisfying
$\left\vert H\right\vert \geq\left(  1-\beta\right)  n$ and one of the conditions:

\emph{(i)} $\mu\left(  H\right)  >\left(  1-1/r\right)  \left\vert
H\right\vert ;$

\emph{(ii)} $\mu\left(  H\right)  >\left(  1-1/r-b\right)  \left\vert
H\right\vert $ and $\delta\left(  H\right)  >\left(  1-1/r-1/r^{4}\right)
\left\vert H\right\vert .$

If condition \emph{(i)} holds, by Theorem \ref{th1} we have
\begin{align*}
js_{r+1}\left(  G\right)   &  \geq js_{r+1}\left(  H\right)  \geq
\frac{\left\vert H\right\vert ^{r-1}}{r^{2r+4}}\geq\left(  1-\beta\right)
^{r-1}\frac{n^{r-1}}{r^{2r+4}}=\left(  1-\frac{4b}{1/r^{4}-b}\right)
^{r-1}\frac{n^{r-1}}{r^{2r+4}}\\
&  >\left(  1-\frac{1}{r^{2}}\right)  ^{r-1}\frac{n^{r-1}}{r^{2r+4}}>\left(
1-\frac{r-1}{r^{2}}\right)  \frac{n^{r-1}}{r^{2r+4}}>\frac{n^{r-1}}{r^{2r+5}},
\end{align*}
implying condition \emph{(a)} and completing the proof.

Suppose now that condition \emph{(ii)} holds. If $H$ contains a $K_{r+1},$ by
Fact \ref{leKd}, we see that
\[
js_{r+1}\left(  G\right)  \geq js_{r+1}\left(  H\right)  \geq\frac{\left\vert
H\right\vert ^{r-1}}{r^{r+3}}\geq\left(  1-\beta\right)  ^{r-1}\frac{n^{r-1}%
}{r^{r+3}}>\frac{n^{r-1}}{2^{r-1}r^{r+3}}>\frac{n^{r-1}}{r^{2r+5}},
\]
implying condition \emph{(a)}.

If $H$ contains no $K_{r+1},$ by Fact \ref{tstab}, $H$ contains an induced
$r$-partite subgraph $H_{0}$ satisfying $\left\vert H_{0}\right\vert >\left(
1-3b^{1/3}\right)  \left\vert H\right\vert $ and $\delta\left(  H_{0}\right)
>\left(  1-6b^{1/3}\right)  \left\vert H\right\vert .$ Now from%
\[
\beta=\frac{4b}{1/r^{4}-b}\leq\frac{4b}{1/r^{4}-1/\left(  2^{10}r^{6}\right)
}\leq8r^{4}b<b^{1/3},
\]
we deduce that%
\[
\left\vert H_{0}\right\vert \geq\left(  1-3b^{1/3}\right)  \left\vert
H\right\vert \geq\left(  1-3b^{1/3}\right)  \left(  1-\beta\right)  n>\left(
1-4b^{1/3}\right)  n
\]
and%
\[
\delta\left(  H_{0}\right)  \geq\left(  1-6b^{1/3}\right)  \left\vert
H\right\vert \geq\left(  1-7b^{1/3}\right)  \left(  1-\beta\right)  n>\left(
1-7b^{1/3}\right)  n.
\]
Thus condition \emph{(b)} holds, completing the proof.
\end{proof}

\bigskip

\begin{proof}
[\textbf{Proof of Theorem \ref{th2.2}}]Let $G$ be a graph of order $n$ with
$\mu\left(  G\right)  >\left(  1-1/r-b\right)  n$. If $G$ contains no
$K_{r+1},$ then condition \emph{(b) }follows from Fact \ref{tstab}; thus we
assume that $G$ contains a $K_{r+1}.$ If
\begin{equation}
\delta\left(  G\right)  >\left(  1-1/r-1/r^{4}\right)  n, \label{mind3}%
\end{equation}
then Fact \ref{thv4} implies condition \emph{(a)}.

Thus, we shall assume that (\ref{mind3}) fails. Then, letting
\begin{equation}
\alpha=1/r^{4}-b,\text{ \ \ }\beta=4b/\alpha,\text{ \ \ }\gamma=1-1/r-b,\text{
\ \ }K=0, \label{co10}%
\end{equation}
we easily see that
\begin{equation}
\beta=\frac{4b}{1/r^{4}-b}\leq\frac{1}{2},\text{ \ \ \ }\delta\left(
G\right)  \leq\left(  \gamma-\alpha\right)  n, \label{co11}%
\end{equation}
and
\begin{equation}
\mu\left(  G\right)  >\left(  1-1/r-b\right)  n=\gamma n. \label{co12}%
\end{equation}
Given (\ref{co10}), (\ref{co11}) and (\ref{co12}), Theorem \ref{thv3} implies
that, for $n\geq r^{20}$, $G$ contains an induced subgraph $H$ satisfying
$\left\vert H\right\vert \geq\left(  1-\beta\right)  n$ and one of the conditions:

\emph{(i)} $\mu\left(  H\right)  >\left(  1-1/r\right)  \left\vert
H\right\vert ;$

\emph{(ii)} $\mu\left(  H\right)  >\left(  1-1/r-b\right)  \left\vert
H\right\vert $ and $\delta\left(  H\right)  >\left(  1-1/r-1/r^{4}\right)
\left\vert H\right\vert .$

If condition \emph{(i)} holds, Theorem \ref{th2} implies that $H$ contains a%
\[
K_{r}^{+}\left(  \left\lfloor 2c\ln\left\vert H\right\vert \right\rfloor
,\ldots,\left\lfloor 2c\ln\left\vert H\right\vert \right\rfloor ,\left\lceil
\left\vert H\right\vert ^{1-2cr^{3}}\right\rceil \right)  .
\]
Now condition \emph{(a) }follows in view of $2c\ln\left\vert H\right\vert
\geq2c\ln\frac{n}{2}>c\ln n$ and%
\[
\left\vert H\right\vert ^{1-2cr^{3}}\geq\left(  \frac{n}{2}\right)
^{1-2cr^{3}}\geq\frac{1}{2}n^{1-2cr^{3}}>n^{1-\sqrt{c}},
\]
completing the proof.

Suppose now that condition \emph{(ii)} holds. If $H$ contains a $K_{r+1},$ by
Fact \ref{thv4}, $H$ contains a%
\[
K_{r}^{+}\left(  \left\lfloor 2c\ln\left\vert H\right\vert \right\rfloor
,\ldots,\left\lfloor 2c\ln\left\vert H\right\vert \right\rfloor ,\left\lceil
\left\vert H\right\vert ^{1-2cr^{3}}\right\rceil \right)  .
\]
This implies condition \emph{(a) }in view of $2c\ln\left\vert H\right\vert
\geq2c\ln\frac{n}{2}>c\ln n$ and%
\[
\left\vert H\right\vert ^{1-2cr^{3}}\geq\left(  \frac{n}{2}\right)
^{1-2cr^{3}}\geq\frac{1}{2}n^{1-2cr^{3}}>n^{1-\sqrt{c}}.
\]

If $H$ contains no $K_{r+1},$ the proof is completed as the proof of Theorem
\ref{th1.2}.
\end{proof}

\subsubsection*{Concluding remarks}

It is not difficult to show that if $G$ is a graph of order $n,$ then the
inequality $e\left(  G\right)  >e\left(  T_{r}\left(  n\right)  \right)  $
implies the inequality $\mu\left(  G\right)  >\mu\left(  T_{r}\left(
n\right)  \right)  .$ Therefore, Theorems \ref{th1}-\ref{th3.2} imply the
corresponding nonspectral extremal results with narrower ranges of the
parameters.\medskip

Finally, a word about the project mentioned in the introduction: in this
project we aim to give wide-range results that can be used further, adding
more integrity to spectral extremal graph theory.

\end{document}